\documentclass[11pt]{article}
\usepackage{amscd,amsmath, amssymb, fancyhdr}

\newcommand{\version}{version 3.0,\ \   Febr. 18, 2021}

\renewcommand{\leftmark}{{\scriptsize LCK manifolds admitting a holomorphic
conformal flow}}

\numberwithin{equation}{section}

\def\eqref#1{(\ref{#1})}

\newcommand{\arrow}{{\:\longrightarrow\:}}
\newcommand{\Z}{{\Bbb Z}}
\newcommand{\C}{{\Bbb C}}
\newcommand{\R}{{\Bbb R}}

\def\1{\sqrt{-1}\:}

\newcommand{\cntrct}                
{\hspace{2pt}\raisebox{1pt}{\text{$\lrcorner$}}\hspace{2pt}}

\renewcommand{\phi}{\varphi}
\renewcommand{\epsilon}{\varepsilon}
\renewcommand{\geq}{\geqslant}

\newcommand{\Aut}{\operatorname{Aut}}

\newcommand{\Diff}{\operatorname{Diff}}

\newcounter{Mycounter}[section]
\newcounter{lemma}[section]
\newcounter{claim}[section]
\newcounter{sublemma}[section]
\newcounter{corollary}[section]
\renewcommand{\thecorollary}{{Corollary \thesection.\arabic{corollary}}}
\newcommand{\corollary}{%
     \setcounter{corollary}{\value{Mycounter}}
     \refstepcounter{corollary}
     \stepcounter{Mycounter}
     {\noindent \bf \thecorollary:\ }}

\newcounter{theorem}[section]
\renewcommand{\thetheorem}{{Theorem \thesection.\arabic{theorem}}}
\newcommand{\theorem}{%
     \setcounter{theorem}{\value{Mycounter}}
     \refstepcounter{theorem}
     \stepcounter{Mycounter}
     {\noindent \bf \thetheorem:\ }}

\newcounter{conjecture}[section]
\newcounter{proposition}[section]
\newcounter{definition}[section]
\renewcommand{\thedefinition}
       {{Definition~\thesection.\arabic{definition}}}
\newcommand{\definition}{%
     \setcounter{definition}{\value{Mycounter}}
     \refstepcounter{definition}
     \stepcounter{Mycounter}
     {\noindent \bf \thedefinition:\ }}

\newcounter{example}[section]
\newcounter{remark}[section]
\renewcommand{\theremark}{{Remark \thesection.\arabic{remark}}}
\newcommand{\remark}{%
     \setcounter{remark}{\value{Mycounter}}
     \refstepcounter{remark}
     \stepcounter{Mycounter}
     {\noindent \bf \theremark:\ }}

\newcounter{problem}[section]
\newcounter{question}[section]
\makeatletter

\@addtoreset{equation}{section}
\@addtoreset{footnote}{section}
\makeatother

\def\blacksquare{\hbox{\vrule width 5pt height 5pt depth 0pt}}
\def\endproof{\hfill\blacksquare}

\begin{document}
\begin{center}
{\LARGE\bf
Locally conformally K\"ahler manifolds\\[3mm]
admitting a holomorphic conformal flow}\\[3mm]

Liviu Ornea\footnote{Partially supported by a PN2-IDEI grant, nr. 525.}
and Misha Verbitsky\footnote{Partially supported by
RFBR grant 10-01-93113-NCNIL-a,
RFBR grant 09-01-00242-a, AG Laboratory SU-HSE, RF government 
grant, ag. 11.G34.31.0023, and Science Foundation of 
the SU-HSE award No. 10-09-0015.\\
Both authors thank Oberwolfach Institute for funding a RIP visit in April
2010 when part of this research was done.


}

\end{center}

{\small
\hspace{0.15\linewidth}
\begin{minipage}[t]{0.7\linewidth}
{\bf Abstract} \\
A manifold $M$ is locally conformally K\"ahler (LCK)
if it admits a  K\"ahler covering $\tilde M$ with
monodromy acting by holomorphic homotheties.  Let $M$
be an LCK manifold admitting a holomorphic conformal flow of
diffeomorphisms, lifted to a non-isometric
homothetic flow on $\tilde M$. We show that
$M$ admits an automorphic potential. This version
is added 10 years after publication to correct
some errors in the original.
\end{minipage}
}

\tableofcontents


\section{Introduction}


\subsection{Conformal automorphisms of LCK manifolds}

Locally conformally K{\"a}hler (LCK)
manifolds are, by definition,
complex manifolds of $\dim_\C>1$
admitting a K{\"a}hler covering $\tilde M$
with deck transformations acting by holomorphic homotheties.
The {\bf monodromy group} of an LCK manifold $M$ is
the deck transform group of the smallest K\"ahler
covering $\tilde M \arrow M$.

This condition is equivalent to the existence of a global {\em closed}
one-form $\theta$ (called {\bf the Lee form}) such that the fundamental
two-form $\omega$ satisfies $d\omega=\theta\wedge\omega$.

We shall always assume
that $M$ is not globally conformally equivalent to
a K\"ahler manifold.

In the present paper, we prove that any
compact LCK manifold which admits a holomorphic conformal flow (\emph{i.e.} an $\R$-action by holomorphic transformations of $\omega$ which are holomorphic with respect to the complex structure $I$ of $M$) 
has monodromy $\Z$, provided that this
flow does not act by isometries on the K\"ahler covering
(\ref{_confo_then_Z_Theorem_}).

An especially interesting class of LCK manifold is
called {\bf LCK manifolds with potential} (see
Subsection \ref{_LCK_pote_Subsection_}). These are
manifolds with the K\"ahler metric on $\tilde M$
admitting a K\"ahler potential which is automorphic
with respect to the action of the monodromy group.

In the present paper, we characterize LCK
manifolds with potential in terms of a
holomorphic conformal flow. We prove that $M$
is an LCK manifold with potential if and only if
it admits a holomorphic conformal flow
which does not act by isometries on the K\"ahler covering
(\ref{_potential_flow_equivalent_}).

\subsection{Vaisman manifolds}

\definition
A {\bf Vaisman manifold} is an LCK manifold $(M, \omega, \theta)$
with $\nabla \theta=0$,
and $\nabla$ the Levi-Civita connection.

\hfill

Compact Vaisman manifolds can be characterized in terms of
their automorphism group.

\hfill

\theorem\label{_Kami_Or_Theorem_}
(\cite{_Kamishima_Ornea_})
Let $(M, \omega)$ be a compact LCK manifold admitting a
holomorphic, conformal action of $\C$ which
lifts to an action by non-trivial homotheties on its
K\"ahler covering. Then $(M, \omega)$ is conformally
equivalent to a Vaisman manifold. \endproof

\hfill

This characterization is superficially similar to
the one given in the present paper for LCK manifolds
with the potential (\ref{_potential_flow_equivalent_}).
However, we ask for a holomorphic conformal
$S^1$-action only, and Kamishima-Ornea theorem postulates
existence of a holomorphic conformal $\C$-action.

\hfill

Vaisman manifolds are especially important 
because their topology is easy to control. As shown in
\cite{_OV:Immersion_}, any Vaisman manifold
is diffeomorphic to a locally trivial
elliptic fibration over a projective orbifold.

\subsection{LCK manifolds with potential}
\label{_LCK_pote_Subsection_}

Let $(\tilde M, \tilde\omega)$ be a K\"ahler covering
of an LCK manifold $M$, and let $\Gamma$ be the deck transform
group of $[\tilde M:M]$.
Denote by $\chi:\; \Gamma \arrow \R^{>0}$
the corresponding character of $\Gamma$,
defined through the scale factor of $\tilde \omega$:
\begin{equation}\label{_character_chi_defi_Equation_}
  \gamma^*\tilde\omega=\chi(\gamma) \tilde\omega, \ \
  \forall \gamma\in\Gamma.
\end{equation}

\definition
A differential form $\alpha$ on $\tilde M$ is called
{\bf automorphic} if $\gamma^*\alpha =
\chi(\gamma)\alpha$, where  $\chi:\; \Gamma \arrow \R^{>0}$
is the character of $\Gamma$ defined above.

\hfill

The K\"ahler form $\tilde\omega$ on every K\"ahler covering
$(\tilde M, \tilde\omega)$
of an LCK manifold is by definition automorphic.

\hfill

\definition\label{_LCK_with_pote_Definition_}
Let $(\tilde M, \tilde\omega)$ be a K\"ahler covering
of an LCK manifold $M$. We say that $M$ is
{\bf an LCK manifold with an automorphic potential}
if $\tilde \omega = dd^c \phi$, for some automorphic
function $\phi$ on $\tilde M$.

\hfill

As shown {\em e.g.} in \cite{_Verbitsky_vanishing_}, a
Vaisman manifold has an automorphic potential,
which can be written down explicitly as
$|\pi^* \theta|^2$,
where $\pi^*\theta$ is the lift of the Lee
form to the K\"ahler covering of $M$,
and $|\cdot|$ the metric associated with its
K\"ahler form.

\hfill

In \cite{_OV:Potential_}, a definition
of {\bf an LCK manifold with potential} was given.
In this definition, in addition to having an
automorphic potential $\phi$, the function
$\phi:\; \tilde M\arrow \R$ was assumed to be
{\bf proper}, that is, with compact fibers.

As shown in \cite[Proposition 2.15, Theorem 2.16]{_OV_Positivity_},
any complex manifold with an LCK metric with
automorphic potential admits another LCK
metric which also has an automorphic
potential, $\phi':\; \tilde M\arrow \R$,
but $\phi'$ is proper and positive. In the published version of the
present paper\footnote{Mathematische Zeitschrift, Volume 273, Issue 3 (2013), Page 605-611}, 
we claimed a stronger result, showing
that the monodromy of any LCK manifold
with potential is $\Z$ (\ref{_confo_then_Z_Theorem_}).
This is now known to be false (\cite[Theorem 3.4]{_OV:rank_}).

\hfill

We showed in \cite{_OV_Top_Potential_}
that any compact LCK manifold
with automorphic  potential can be obtained as a
deformation of a Vaisman manifold.  Many of the known examples of LCK
manifolds
are Vaisman (see \cite{belgun} for a complete list of Vaisman compact
complex surfaces), but there are also non-Vaisman ones: one of the Inoue
surfaces (see \cite{belgun}, \cite{_Tricerri_}), its higher-dimensional
generalization in \cite{_Oeljeklaus_Toma_}, the non-diagonal Hopf
manifolds $(\C^{N}\setminus\{0\})/\langle A\rangle$ with $A$ linear, with
eigenvalues smaller than $1$ in absolute value (see \cite{GO},
\cite{_OV:Potential_}),
and the new examples found in \cite{FP} on
parabolic and hyperbolic Inoue surfaces.

\hfill

Compact LCK manifolds with potential are embeddable in
Hopf manifolds, see \cite{_ov:MN_}. The existence
of an automorphic potential leads to important topological
restrictions on the fundamental group, see
\cite{_OV_Top_Potential_} and \cite{_Kokarev_Kotschick:Fibrations_}.

The class of compact complex manifolds admitting
an LCK metric with potential is
stable under small complex deformation
(\cite[Theorem 2.6]{_OV:Potential_}).
This statement should be considered as an LCK
analogue of Kodaira's K\"ahler
stability theorem.
It was the original way to construct LCK metrics
on some non-Vaisman manifolds, such as the
Hopf manifolds not admitting a Vaisman structure
(\cite{GO}).

\hfill

In \cite{_OV_Automorphisms_}, it was shown that
LCK manifolds with automorphic
potential can be characterized in terms of existence of a
particular subgroup of automorphisms. To state the result, we need to
introduce the {\em weight bundle} $L\rightarrow M$
associed to the representation $\mathrm{GL}(2n,\R)\ni
A\mapsto |\det A|^{1/n}$. It is endowed with the flat connection form
$-\frac 12 \theta$, thus producing a local system.  The holonomy
of this local system is precisely the
monodromy group of $M$.

\hfill

\theorem\label{_S^1_main_Theorem_}
(\cite[Theorem 1.8]{_OV_Automorphisms_})
Let $M$ be a compact complex manifold, equipped
with a holomorphic $S^1$-action and a LCK metric
(not necessarily compatible). Suppose that
the weight bundle $L$, restricted to a general
orbit of this $S^1$-action, is non-trivial
as a $1$-dimensional local system. Then $M$ admits
a LCK metric with an automorphic potential.

\hfill

\remark
The converse statement is true as well
(\ref{_potential_flow_equivalent_}).
In the present paper we prove that an
LCK manifold $M$ with an automorphic
potential always admits a holomorphic, conformal
$S^1$-action which lifts to an action by
non-trivial homotheties on its covering.

\hfill

As shown in \cite[Corollary 1.11]{_OV_Top_Potential_},
\ref{_S^1_main_Theorem_} implies the following corollary.

\hfill

\corollary
Let $M$ be a compact LCK manifold  of complex dimension
$n\geq 3$. Suppose that
the weight bundle $L$ restricted to a general
orbit of this $S^1$-action is non-trivial
as a 1-dimensional local system. Then $M$
is diffeomorphic to a Vaisman manifold, and admits
a holomorphic embedding to a Hopf manifold.


\section{Holomorphic conformal flows on LCK manifolds}


\subsection{Holomorphic conformal flow and monodromy}

Let $M$ be an LCK manifold, $\tilde M$ its minimal
K\"ahler covering, and $\chi:\; \Gamma \arrow \R^{>0}$ the
character defined through the scale factor as in 
\eqref{_character_chi_defi_Equation_}.
Observe first that because we work with the minimal covering, the
character $\chi$ is injective and hence the monodromy $\Gamma$ can be
viewed as a  subgroup of the
multiplicative group $\R^{>0}$. As such, the monodromy group is
abelian and torsion-free.

\hfill

\theorem\label{_confo_then_Z_Theorem_}
Let $M$ be a compact LCK manifold, and
$\rho:\; \R \times M \arrow M$ a holomorphic conformal flow
of diffeomorphisms on $M$. Assume that $\rho$
is lifted to a flow of non-isometric homotheties
on the smallest K\"ahler covering $\tilde M$ of $M$.
Let $\Gamma\subset \Aut M$ be the deck transform
group (that is, the monodromy), $G$ the
closure of the group generated by $\rho$
in $\Diff(M)$, and $\tilde G$ its preimage in
$\Diff(\tilde M)$. Then $\tilde G$ is connected
and contains $\Gamma$.

\hfill

{\bf Proof:} 
From the exact sequence
\[
0 \arrow \Gamma \arrow \tilde G \arrow G \arrow 0
\]
it is clear that $\Gamma\subset \tilde G$. It
remains only to prove that $\tilde G$ is connected.

Notice that any conformal holomorphic
map $\phi:\; V \arrow W$ of K\"ahler manifolds of dimension
$>1$ is a homothety. Indeed,
the pullback $\phi^* \omega_W$ of the K\"ahler form on $W$ under
a holomorphic morphism is closed. Since $\phi^*\omega_W = f \omega_V$,
for some positive function $f$ on $V$, this gives 
$df\wedge \omega_V=0$, hence $df=0$. 
Then $\tilde G$ is a group of holomorphic
homotheties (it was essential to know that the flow is
holomorphic to be able to deduce that its lift, formed by
conformalities with respect to the K\"ahler metric,
contains in fact only homotheties). 

Observe that every connected component of $\tilde G$ contains an isometry.
Indeed, take an element $\tilde a\in \tilde G$
(where $\tilde a$ is a lift of an $a\in G$). It acts on
the K\"ahler form $\tilde \omega$ as $\tilde a^*\tilde\omega =C_a\cdot
\tilde\omega$,
$C_a=\text{const.}$ Consider the element $c:=\rho_{C_a^{-1}}$ of the
conformal flow which satisfies $\tilde c^*\tilde\omega=C_a^{-1}\tilde\omega$.
Let $\tilde b=\tilde\rho_{C_a^{-1}}\circ\tilde a$. Then
$\tilde b$ is an isometry with respect to $\tilde\omega$ in the
same component with $\tilde a$. 

Let $\tilde G_0\subset \tilde G$ be the group of isometries,
and $G_0$ its image in $G$.
Since $\tilde M$ is the minimal K\"ahler covering, 
the natural projection $\tilde G_0 \arrow G_0$ is an isomorphism
(otherwise, we could take a quotient by the kernel of this projection,
and obtain a smaller K\"ahler covering). 
Therefore, each $g\in G_0$ admits a unique lifting to $\tilde G$.
Since each connected component of $\tilde G$ surjectively
projects to $G$, each of these components would contain
a preimage of $G_0$. This is impossible, because 
each  $g\in G_0$ has a unique preimage. Therefore,
$\tilde G$ is connected.

%
%
%
%

\endproof

\subsection{Holomophic conformal flow and automorphic
potential}

 Applying \cite{_OV_Automorphisms_}, we
obtain that a compact LCK manifold which satisfies the
assumptions of \ref{_confo_then_Z_Theorem_}
always admits an automorphic potential.
This gives the following corollary:

\hfill

\corollary\label{_automo_pote_from_confo_Corollary_}
Let $M$ be a compact LCK manifold, and
$\rho:\; \R \times M \arrow M$ a holomorphic conformal flow
of diffeomorphisms on $M$. Assume that $\rho$
is lifted to a flow of non-isometric homotheties
on the K\"ahler covering $\tilde M$ of $M$.
Then $M$ admits an authomorphic potential.

\hfill

{\bf Proof} Let $\omega_0$ be the Gauduchon Hermitian form associated with
$\omega$. It is a Hermitian form which is conformally
equivalent to $\omega$ and satisfies 
$dI d(\omega^{n-1})=0$, where $n = \dim_\C M$.
As shown in \cite{_Gauduchon_1984_}, such a form
always exists, and is unique up to a constant. The constant
can be chosen in such a way that
$\int_M \omega^n_0=1$, and
in this case $\omega_0$ is preserved
by any conformal holomorphic diffeomorphism.

This implies that $\omega_0$ is $\rho$-invariant.
Let $G$ be the closure of $\rho(\R)$ in the group
of all holomorphic isometris of $M$. Since this
group is a compact Lie group (see {\em e.g.}  
\cite{_Myers_Steenrod_}), $G$ is also a compact
Lie group, which is obviously commutative,
hence isomorphic to a torus. Choosing an 
appropriate 1-dimensional subgroup in $G$,
we arrive in the situation described by 
\cite{_OV_Automorphisms_}: $M$ is equipped
with a holomorphic action of a circle, which
lifts to non-isometric homotheties of
its K\"ahler covering. From  \cite[Theorem
1.8 ]{_OV_Automorphisms_}, we now obtain that 
$M$ admits an automorphic potential.
\endproof

\hfill

\ref{_automo_pote_from_confo_Corollary_} can be used to give
a characterization of LCK manifolds with automorphic potential.

\hfill

\theorem\label{_potential_flow_equivalent_}
Let $(M, I, \omega)$ be a
compact LCK manifold. Then the following
assertions are equivalent.
\begin{description}
\item[(i)] $M$ admits an automorphic potential.
\item[(ii)] The complex manifold
$(M,I)$ admits an LCK metric $\omega'$ with same monodromy,
and a conformal flow of holomorphic diffeomorphisms of $(M,I,\omega')$, which
is lifted to a flow of non-isometric homotheties
of the K\"ahler covering $(\tilde M, \tilde \omega')$ of $(M,\omega')$.
\end{description}

{\bf Proof:} The implication (ii) $\Rightarrow$ (i)
immediately follows from \ref{_automo_pote_from_confo_Corollary_}.  We now show how {(i)}
implies {(ii)}.

Embed $M$ in a Hopf manifold
$H=(\C^N\setminus\{0\})/\langle A\rangle$, where $A$ is a
linear (not necessarily diagonal) operator with
eigenvalues strictly smaller than $1$ in absolute
value. Such an embedding was constructed in
\cite{_OV:Potential_}.
Then $A$ preserves the K\"ahler
covering $\tilde M$ and hence can be considered
as an element of the deck group $\Gamma$. As such it acts
as a homothety on the K\"ahler metric and we can suppose
it is a contraction (otherwise we work with
$A^{-1}$). What we want is to construct (out of $A$)  a
holomorphic flow preserving $\tilde M$.

Recall from \cite{_OV:Potential_} that the metric
completion $\tilde M_c$ of $\tilde M$ is obtained by adding
only one point $z$ (here the existence
of the global potential is crucial). Then $A$ acts
trivially on $z$ and we may consider the local ring
$\mathcal{O}_{\tilde M_c}$ at $z$.

Observe that  $A$ induces an automorphism of the ring
$\mathcal{O}_{\tilde M_c}$, denoted equally by $A$. Then one
easily sees that the formal logarithm of $A$,  $\log A$,
is a derivation of $\mathcal{O}_{\tilde M_c}$ (this follows,
{\em e.g.}, from \cite[p. 209]{_bourbaki_}; it is enough
to show that formally $e^{\log A}=1$). This means that
$\log A$ induces a vector field on $\tilde M$ with
associated flow $e^{t\log A}$. Note that $\log A$ is a holomorphic object
because, as all eigenvalues of $A$ are smaller than $1$ in absolute value,
the corresponding formal series converges. As $M=\tilde M/\langle
A\rangle$, we see that $e^{t\log A}$ projects on a
one-parameter flow on $M$. But, as for $t=1$ the flow on
$\tilde M$ is $A$ which acts trivially on $M$, the orbits
of the projected flow are closed, and hence the projected
flow corresponds to an $S^1$-action on $M$. This action is
holomorphic because $A$ acts holomorphically on $\tilde
M$.

Apply now the averaging (on $S^1$) procedure described in
\cite[2.1]{_OV_Automorphisms_} to obtain a new LCK metric
$\omega'$ on $M$ with respect to which this $S^1$ acts by
holomorphic isometries. We note that the averaging steps
performed do not change the cohomology class of the Lee
form, and hence the new LCK structure has the same
monodromy.

It remains to justify why the lift of this isometric and holomorphic  $S^1$
to $(\tilde M, \tilde \omega')$ is by non-trivial
homotheties. This is because the lifted flow contains $A$
which is a contraction with respect to a certain metric
and hence cannot be an isometry with respect to any
metric.
\endproof

\hfill

\remark
The averaging construction
used in the proof of \ref{_potential_flow_equivalent_}
preserves the class of LCK metrics with potential.
This means that $\omega'$ of
\ref{_potential_flow_equivalent_} (ii) has an automorphic
potential if and only if $\omega$ has one.

\hfill

%
%

\noindent {\bf Acknowledgment.} The authors are grateful
to the referee for carefully reading the paper and for the referee's
very useful suggestions.

{\small

\noindent {\sc Liviu Ornea\\
University of Bucharest, Faculty of Mathematics, \\14
Academiei str., 70109 Bucharest, Romania. \emph{and}\\
Institute of Mathematics ``Simion Stoilow" of the Romanian Academy,\\
21, Calea Grivitei Street
010702-Bucharest, Romania }\\
\tt Liviu.Ornea@imar.ro, \ \ lornea@gta.math.unibuc.ro

\hfill

\noindent {\sc Misha Verbitsky\\
Laboratory of Algebraic Geometry, \\
Faculty of Mathematics, NRU HSE,\\
7 Vavilova Str. Moscow, Russia
}
}


{\scriptsize
\begin{thebibliography}{100}



\bibitem[B]{belgun} F.A. Belgun, {\em On the metric structure of
non-K{\"a}hler complex surfaces}, Math. Ann. {\bf 317} (2000),
1--40.



\bibitem[Bo]{_bourbaki_} Bourbaki, Groupes et alg\`ebres de Lie, Ch. 2-3,
Hermann,
Paris, 1972.

\bibitem[DO]{drag}
S.  Dragomir and L.  Ornea,  Locally conformal
K{\"a}hler
geometry, Progress in Math. {\bf 155},   Birkh{\"a}user, Boston, Basel, 1998.


\bibitem[FP]{FP} A. Fujiki and M. Pontecorvo, \emph{Anti-self-dual
bihermitian structures on Inoue surfaces}, arXiv:0903.1320.

\bibitem[Ga]{_Gauduchon_1984_} 
Gauduchon, 
P. {\em La 1-forme de torsion d'une variete 
hermitienne compacte}, Math. Ann., 267 (1984), 495-518.


\bibitem[GO]{GO} 
P. Gauduchon and L. Ornea, \emph{Locally conformal
K\"ahler metrics on Hopf surfaces}, Ann. Inst. Fourier, {\bf 48} (1998),
1107--1127.


\bibitem[KO]{_Kamishima_Ornea_}
Y. Kamishima, L. Ornea,
{\em Geometric flow on compact locally conformally Kahler manifolds},
Tohoku Math. J.  {\bf 57}  (2005),  no. 2, 201--222.


\bibitem[KK]{_Kokarev_Kotschick:Fibrations_}
G. Kokarev, D. Kotschick, \emph{Fibrations and
fundamental
groups of K\"ahler-Weyl manifolds},  Proc. Amer. Math. Soc.  {\bf 138} 
(2010),   997--1010. arXiv:0811.1952.


\bibitem[MS]{_Myers_Steenrod_} S.B. Myers, N. Steenrod, \emph{The
group of isometries of Riemannian manifolds}, Ann. Math. {\bf 40}
(1939), 400--416.



\bibitem[OT]{_Oeljeklaus_Toma_}
K. Oeljeklaus, M. Toma, \emph{Non-K\"ahler compact
complex manifolds associated to number fields}, Ann. Inst. Fourier {\bf
55} (2005), 1291--1300.


\bibitem[OV1]{_OV:Structure_} L. Ornea and M. Verbitsky,
  \emph{Structure theorem for compact Vaisman manifolds},
  Math. Res. Lett., {\bf 10} (2003), 799--805.


\bibitem[OV2]{_OV:Immersion_}
L. Ornea and M. Verbitsky,
\emph{An immersion
theorem for
compact Vaisman manifolds},  Math. Ann. {\bf 332}  (2005),
121--143. math.AG/0306077


\bibitem[OV3]{_OV:Potential_}
L. Ornea and M. Verbitsky, \emph{Locally
conformal K\"ahler manifolds with potential}, Math. Ann. {\bf 348} (2010), 25--33. 
arXiv:math/0407231


\bibitem[OV4]{_ov:MN_}
L. Ornea and M. Verbitsky, {\em
 Morse-Novikov cohomology of locally conformally K\"ahler  manifolds},
J.Geom.Phys. {\bf 59},(2009), 295--305. arXiv:0712.0107


\bibitem[OV5]{_OV_Top_Potential_}
L. Ornea and M. Verbitsky, {\em Topology of locally conformally K\"ahler
manifolds with potential}, Int. Math. Res. Not. IMRN, {\bf 4} (2010),
117--126. arXiv:0904.3362.

\bibitem[OV6]{_OV_Automorphisms_}
L. Ornea and M. Verbitsky,
{\em Automorphisms of locally conformally K\"ahler
manifolds with potential}, Int. Math. Res. Not, IMRN, doi:10.1093/imrn/rnr053. arXiv:0906.2836.

\bibitem[OV7]{_OV:rank_}
L. Ornea and M. Verbitsky,
{\em LCK rank of locally conformally Kahler manifolds with potential},
J. Geom. Phys. 107 (2016), 92-98.

\bibitem[OV8]{_OV_Positivity_}
L. Ornea and M. Verbitsky,  {\em 
Positivity of LCK potential},
J. Geom. Anal. 29 (2019), no. 2, 1479-1489. 


\bibitem[Tr]{_Tricerri_} F. Tricerri, {\it Some examples of locally
conformal K{\"a}hler manifolds}, Rend.  Sem.  Mat.  Univ.
Politec.  Torino {\bf 40} (1982), 81--92.



\bibitem[Ve1]{_Verbitsky_vanishing_}
M. Verbitsky, {\em Theorems on the
vanishing of cohomology for
locally conformally hyper-K\"ahler manifolds}, Proc. Steklov Inst. Math.
{\bf 246} (2004) 54--78. arXiv:math/0302219.



\end{thebibliography}
}
\end{document}